January 12, 2017

# POWER SERIES WITH SKEW-HARMONIC NUMBERS, DILOGRITHMS, AND DOUBLE INTEGRALS

Khristo N. Boyadzhiev

Department of Mathematics and Statistics,
Ohio Northern University,
Ada, OH 45810, USA
k-boyadzhiev@onu.edu

**Abstract.** The skew-harmonic numbers are the partial sums of the alternating harmonic series, i.e. the expansion of log2. We evaluate in closed form various power series and numerical series with skew-harmonic numbers. This provides a simultaneous solution of two recent problems by Ovidiu Furdui in the American Mathematical Monthly and the College Mathematics Journal. We also present and discuss representations involving the dilogarithm and the trilogaithm which are related to our results. Finally, we provide the evaluations of several double integrals in terms of classical constants.

**1. Introduction.**

We first define the skew-harmonic numbers

$$H_n^- = 1 - \frac{1}{2} + \ldots + \frac{(-1)^{n-1}}{n},$$

for $n = 1, 2, \ldots,$. For convenience, let $H_0^- = 0$. The skew-harmonic numbers are the partial sums in the expansion

$$\log 2 = \sum_{n=1}^{\infty} \frac{(-1)^{n-1}}{n}.$$

In this paper we consider various power series (generating functions) with skew-harmonic numbers in the coefficients. These series are evaluated in terms of logarithms, dilogarithms and trilogarithms. All this is presented in Section 2. In this process we obtain interesting relations for the dilogaithm and the trilogithm which are discussed in Section 3. In particular, in section 3 we give the Taylor expansions in the variable $x$ for the functions ($\mu$ - a parameter)

$$\text{Li}_m\left(\frac{(1+\mu)x}{1+x}\right) - \text{Li}_m\left(\frac{x}{1+x}\right), \quad m = 2,3,$$

where the Taylor coefficients are written in terms of skew-harmonic numbers.

In section 4 we connect our result to the values of certain double integrals.

This work was inspired by three problems posed recently by Ovidiu Furdui in the Fibonacci Quarterly [4], the American Mathematical Monthly [6], and the College Mathematics Journal [7]; see below (15) and (16). These problems and several similar to them have also appeared in the recent book [8].

One of our results is computing the generating function

$$g(z) = \sum_{n=0}^{\infty} \left(H_n^- - \log 2\right)^2 z^n$$

($|z| \leq 1$) in closed form, which extends Furdui's cases (15) and (16).

The skew-harmonic numbers $H_n^-$ are related to the digamma function $\psi(x) = \dfrac{d}{dx}\log\Gamma(x)$ by the equation

$$\psi\left(\frac{n+1}{2}\right) - \psi\left(\frac{n}{2}\right) = 2(-1)^{n-1}\left(\log 2 - H_{n-1}^-\right) \tag{1}$$

which shows the relationship between our work and De Doelder's results in section 5 of his paper [5].

The skew-harmonic numbers can also be represented by a combination of harmonic numbers. It is easy to check that for every $n = 1, 2, \ldots$,

$$H_{2n}^- = H_{2n} - H_n, \quad H_{2n+1}^- = H_{2n+1} - H_n.$$

## 2. Power series with skew-harmonic numbers

We need a simple lemma about power series.

**Lemma 1.** *For any power series $f(t) = a_0 + a_1 t + a_2 t^2 + \ldots$ and for small enough parameter $\lambda$*

$$\frac{f(t)}{1-\lambda t} = \sum_{n=0}^{\infty} \left\{ \sum_{k=0}^{n} \lambda^{n-k} a_k \right\} t^n .$$

The proof follows from Cauchy's rule for multiplication of power series.

Applying the lemma with $\lambda = 1$ to the series

$$\log(1+t) = \sum_{n=1}^{\infty} \frac{(-1)^{n-1}}{n} t^n ,$$

yields immediately the generating function for $H_n^-$,

$$\frac{\log(1+t)}{1-t} = \sum_{n=1}^{\infty} H_n^- t^n \qquad (2)$$

Our purpose is to construct various power series with $H_n^-$ in the coefficients. For instance, combining (2) with

$$\frac{\log 2}{1-t} = \sum_{n=0}^{\infty} (\log 2) t^n ,$$

we find the interesting representation

$$\sum_{n=0}^{\infty} \left( H_n^- - \log 2 \right) t^n = \frac{1}{1-t} \log\left( \frac{1+t}{2} \right) , \qquad (3)$$

where the series converges for $|t| < 1$ and even for $t = 1$, namely,

$$\sum_{n=0}^{\infty} \left( H_n^- - \log 2 \right) = \frac{-1}{2}. \qquad (4)$$

This series is a convergent alternating series where $|H_n^- - \log 2| \leq \dfrac{1}{n+1}$,

and $\lim\limits_{t \to 1} \dfrac{1}{1-t} \log\left( \dfrac{1+t}{2} \right) = \dfrac{-1}{2}$.

In the next proposition we shall use the dilogarithm function as defined in [11],

$$\operatorname{Li}_2(t) = \sum_{n=1}^{\infty} \frac{t^n}{n^2}, \quad |t| \leq 1.$$

**Proposition 2.** *For all* $|t| \leq 1, t \neq 1$,

$$\sum_{n=1}^{\infty} H_n^- \frac{t^n}{n} = \operatorname{Li}_2\left(\frac{1-t}{2}\right) - \operatorname{Li}_2\left(\frac{1}{2}\right) - \operatorname{Li}_2(-t) - \log(1-t)\log 2. \tag{5}$$

(This representation appears in a modified form in Hansen table [10, (5.5.27)]; it also follows from formula A.2.8 (1) in [11]).

*Proof.* Notice that

$$\frac{d}{dx}\operatorname{Li}_2\left(\frac{1-x}{2}\right) = \frac{1}{1-x}\log\left(\frac{1+x}{2}\right) = \frac{\log(1+x)}{1-x} - \frac{\log 2}{1-x}. \tag{6}$$

Therefore,

$$\sum_{n=1}^{\infty} H_n^- x^n = \frac{\log(1+x)}{1-x} = \frac{d}{dx}\operatorname{Li}_2\left(\frac{1-x}{2}\right) + \frac{\log 2}{1-x}.$$

Integrating this we find

$$\sum_{n=1}^{\infty} \frac{x^{n+1}}{n+1} H_n^- = \operatorname{Li}_2\left(\frac{1-x}{2}\right) - \operatorname{Li}_2\left(\frac{1}{2}\right) - \log 2 \log(1-x) \tag{7}$$

which reduces to (5), as

$$\sum_{n=1}^{\infty} \frac{x^{n+1}}{n+1} H_n^- = \sum_{n=1}^{\infty} \frac{x^n}{n} H_n^- + \operatorname{Li}_2(-x). \qquad \square$$

**Remark 3**. Equation (5) can be written also in the form

$$\sum_{n=1}^{\infty} \left(H_n^- - \log 2\right) \frac{t^n}{n} = \operatorname{Li}_2\left(\frac{1-t}{2}\right) - \operatorname{Li}_2\left(\frac{1}{2}\right) - \operatorname{Li}_2(-t) \tag{8}$$

by moving the term $\log(1-t)\log 2$ to the LHS and expanding $\log(1-t)$.

Note that (8) is convergent on the entire closed disk $|t| \leq 1$, while (5) is not convergent for $t=1$. Thus

$$\sum_{n=1}^{\infty}\left(H_n^- - \log 2\right)\frac{1}{n} = -\text{Li}_2\left(\frac{1}{2}\right) - \text{Li}_2(-1) = \frac{1}{2}\log^2 2$$

as $\text{Li}_2(-1) = \frac{-\pi^2}{12}$ and $\text{Li}_2\left(\frac{1}{2}\right) = \frac{\pi^2}{12} - \frac{1}{2}\log^2 2$ (cf. Problem 3.14 in [8]).

With $t = -1$,

$$\sum_{n=1}^{\infty}\left(H_n^- - \log 2\right)\frac{(-1)^{n-1}}{n} = \frac{\pi^2}{12} - \frac{\log^2 2}{2} \, . \tag{9}$$

This is equivalent to

$$\sum_{n=1}^{\infty} H_n^- \frac{(-1)^{n-1}}{n} = \frac{\pi^2}{12} + \frac{\log^2 2}{2}, \tag{10}$$

which results also from (5) with $t = -1$.

**Remark 4.** From (3) and (6),

$$\sum_{n=0}^{\infty}\left(H_n^- - \log 2\right) t^n = \frac{d}{dt}\text{Li}_2\left(\frac{1-t}{2}\right).$$

Integrating this equation and adjusting the constant of integration we find

$$\sum_{n=0}^{\infty}\left(H_n^- - \log 2\right)\frac{t^{n+1}}{n+1} = \text{Li}_2\left(\frac{1-t}{2}\right) - \text{Li}_2\left(\frac{1}{2}\right). \tag{11}$$

**Proposition 5.** *For $|t| < 1$,*

$$\sum_{n=1}^{\infty}\left(H_n^-\right)^2 t^n = \frac{1}{1-t}\left\{\text{Li}_2(t) + 2\log 2 \log(1+t) + 2\text{Li}_2\left(\frac{1}{2}\right) - 2\text{Li}_2\left(\frac{1+t}{2}\right)\right\} \tag{12}$$

*and for $|t| \leq 1$ (using the limit of the RHS for $t \to 1$),*

$$\sum_{n=0}^{\infty}\left(H_n^- - \log 2\right)^2 t^n = \frac{1}{1-t}\left\{\text{Li}_2(t) + \log^2 2 - 2\left[\text{Li}_2\left(\frac{1+t}{2}\right) - \text{Li}_2\left(\frac{1}{2}\right)\right]\right\} . \tag{13}$$

For the proof we need the following lemma.

**Lemma 6**. *For every $n = 1, 2, \ldots$, and with $H_n^{(2)} = 1 + \frac{1}{2^2} + \ldots + \frac{1}{n^2}$ we have*

$$2\sum_{k=1}^{n} \frac{(-1)^{k-1}}{k} H_k^- = \left(H_n^-\right)^2 + H_n^{(2)}, \tag{14}$$

*Proof by induction*. For $n = 1$ the equality is clear. If it holds for some $n$, then

$$2\sum_{k=1}^{n+1} \frac{(-1)^{k-1} H_k^-}{k} = 2\sum_{k=1}^{n} \frac{(-1)^{k-1} H_k^-}{k} + 2\frac{(-1)^n H_{n+1}^-}{n+1}$$

$$= \left(H_n^-\right)^2 + H_n^{(2)} + 2\frac{(-1)^n H_{n+1}^-}{n+1}$$

$$= \left(H_{n+1}^- - \frac{(-1)^n}{n+1}\right)^2 + H_{n+1}^{(2)} - \frac{1}{(n+1)^2} + 2\frac{(-1)^n H_{n+1}^-}{n+1}$$

$$= \left(H_{n+1}^-\right)^2 + H_{n+1}^{(2)}. \qquad \square$$

Next we turn to the proof of the proposition. Replacing $t$ by $-t$ in (5) we write

$$\sum_{n=1}^{\infty} \frac{(-1)^n H_n^-}{n} t^n = \text{Li}_2\left(\frac{1+t}{2}\right) - \text{Li}_2\left(\frac{1}{2}\right) - \text{Li}_2(t) - \log(1+t)\log 2,$$

and with the help of lemma 1 we derive

$$\sum_{n=1}^{\infty} \left\{2\sum_{k=1}^{n} \frac{(-1)^{k-1} H_k^-}{k}\right\} t^n = \frac{2}{1-t}\left\{\text{Li}_2(t) + \log(1+t)\log 2 - \text{Li}_2\left(\frac{1+t}{2}\right) + \text{Li}_2\left(\frac{1}{2}\right)\right\}.$$

At the same time,

$$\sum_{n=1}^{\infty} H_n^{(2)} t^n = \frac{\text{Li}_2(t)}{1-t},$$

and from Lemma 6 we obtain (12). After that we write

$$\left(H_n^- - \log 2\right)^2 = \left(H_n^-\right)^2 - (2\log 2) H_n^- + \log^2 2,$$

and correspondingly, with summation starting from $n = 0$,

$$\sum_{n=0}^{\infty} \left(H_n^- - \log 2\right)^2 t^n = \sum_{n=0}^{\infty} \left(H_n^-\right)^2 t^n - 2\log 2 \sum_{n=0}^{\infty} H_n^- t^n + \sum_{n=0}^{\infty} (\log^2 2) t^n.$$

Now (13) follows from (12) and (2).  □

**Remark 7.** Setting $t = -1$ in (13) we compute

$$\sum_{n=0}^{\infty}(-1)^n \left(H_n^- - \log 2\right)^2 = \frac{1}{2}\left\{\text{Li}_2(-1) + \log^2 2 + 2\text{Li}_2\left(\frac{1}{2}\right)\right\} = \frac{\pi^2}{24}. \tag{15}$$

This provides an independent solution to problem 11682 in the American Mathematical Monthly [6]. The series (13) converges also for $t = 1$, as we have the estimate

$$\left(H_n^- - \log 2\right)^2 \leq \frac{1}{(n+1)^2},$$

and computing the limit of the RHS in (13) when $t \to 1$ we find the remarkable result:

$$\sum_{n=0}^{\infty}\left(H_n^- - \log 2\right)^2 = \log 2. \tag{16}$$

This is an independent solution to problem 997 in the College Mathematics Journal [7].

**Proposition 8.** *For all $|x| \leq 1$,*

$$\sum_{n=0}^{\infty}\left(H_n^- - \log 2\right)^2 \frac{x^{n+1}}{n+1} = \int_0^x \frac{\text{Li}_2(t)}{1-t}dt + \log(1-x)\left[2\text{Li}_2\left(\frac{1+x}{2}\right) - \frac{\pi^2}{6}\right] \tag{17}$$

$$+ 2\log(1+x)\left(\log^2(1-x) - \log^2 2\right) + 2\log 2\left[\text{Li}_2\left(\frac{1+x}{2}\right) - \text{Li}_2\left(\frac{1}{2}\right)\right]$$

$$-2\log 2 \log^2(1-x) + 4\log(1-x)\text{Li}_2\left(\frac{1-x}{2}\right) - 4\left[\text{Li}_3\left(\frac{1-x}{2}\right) - \text{Li}_3\left(\frac{1}{2}\right)\right].$$

*In particular, for $x = \pm 1$ we have*:

$$\sum_{n=0}^{\infty}\left(H_n^- - \log 2\right)^2 \frac{1}{n+1} = \frac{3}{2}\zeta(3) - \frac{\pi^2}{6}\log 2 - \frac{1}{3}\log^3 2, \tag{18}$$

$$\sum_{n=0}^{\infty}\left(H_n^- - \log 2\right)^2 \frac{(-1)^{n+1}}{n+1} = \frac{\pi^2}{12}\log 2 - \frac{3}{4}\zeta(3) - \frac{1}{3}\log^3 2. \tag{19}$$

*Proof.* To obtain (17) we integrate for $t$ in (13) from $0$ to $x$. The integral

$$\int_0^x \left[\text{Li}_2\left(\frac{1+t}{2}\right) - \text{Li}_2\left(\frac{1}{2}\right)\right]\frac{dt}{1-t}$$

is solved by parts and the resulting simpler integrals are left to the reader. The integral $\int_0^x \frac{\text{Li}_2(t)}{1-t} dt$ is solved by the techniques in [11] and we find two versions:

$$\int_0^x \frac{\text{Li}_2(t)}{1-t} dt = -2\text{Li}_3\left(\frac{-x}{1-x}\right) - 2\text{Li}_3(x) + \log(1-x)\text{Li}_2(x) + \frac{1}{3}\log^3(1-x)$$

which holds for $-1 \leq x \leq \frac{1}{2}$, and

$$\int_0^x \frac{\text{Li}_2(t)}{1-t} dt = 2\left[\text{Li}_3(1-x) - \zeta(3)\right] - \log(1-x)\left[\text{Li}_2(1-x) + \frac{\pi^2}{6}\right]$$

for $0 \leq x < 1$. The integral is divergent at $t = 1$, but the entire expression on the RHS in the proposition has a limit when $x \to 1$ and we shall compute this limit. First, however, we shall evaluate the series for $x = -1$. Thus

$$\int_0^{-1} \frac{\text{Li}_2(t)}{1-t} dt = \frac{\pi^2}{12} \log 2 - \frac{1}{4} \zeta(3) .$$

Computing the limit for $x \to -1$ of the remaining part in the RHS in (17) we find

$$\sum_{n=0}^{\infty} \left(H_n^- - \log 2\right)^2 \frac{(-1)^{n+1}}{n+1} = \frac{\pi^2}{12} \log 2 - \frac{3}{4} \zeta(3) - \frac{1}{3} \log^3 2 .$$

Here we have used the values

$$\text{Li}_3\left(\frac{1}{2}\right) = \frac{7}{8}\zeta(3) - \frac{\pi^2 \log 2}{12} + \frac{(\log 2)^3}{6} \quad \text{and} \quad \text{Li}_3(-1) = \frac{-3}{4}\zeta(3) .$$

In order to evaluate the series for $x = 1$ we first write

$$\int_0^x \frac{\text{Li}_2(t)}{1-t} dt + \log(1-x)\left[2\text{Li}_2\left(\frac{1+x}{2}\right) - \frac{\pi^2}{6}\right]$$

$$= 2\left[\text{Li}_3(1-x) - \zeta(3)\right] - \log(1-x)\text{Li}_2(1-x) + \log(1-x)\left[2\text{Li}_2\left(\frac{1+x}{2}\right) - \frac{\pi^2}{3}\right],$$

and when $x \to 1$ the limit of this function equals $-2\zeta(3)$.

The remaining part of the RHS in the proposition we rearrange in the form

$$2\log^2(1-x)[\log(1+x)-\log 2]+2\log 2\left[\operatorname{Li}_2\left(\frac{1+x}{2}\right)-\operatorname{Li}_2\left(\frac{1}{2}\right)-\log 2\log(1+x)\right]$$

$$+4\log(1-x)\operatorname{Li}_2\left(\frac{1-x}{2}\right)-4\left[\operatorname{Li}_3\left(\frac{1-x}{2}\right)-\operatorname{Li}_3\left(\frac{1}{2}\right)\right].$$

Evaluating the limit of this function for $x \to 1$ we find its value to be

$$\frac{7}{2}\zeta(3)-\frac{\pi^2}{6}\log 2-\frac{1}{3}\log^3 2 \ .$$

Adding $-2\zeta(3)$ to this we come to the desired result

$$\sum_{n=0}^{\infty}\left(H_n^- - \log 2\right)^2 \frac{1}{n+1} = \frac{3}{2}\zeta(3)-\frac{\pi^2}{6}\log 2-\frac{1}{3}\log^3 2 \ . \qquad \square$$

**Proposition 9.** *For all* $\dfrac{-1}{3} \le x \le 1$,

$$\sum_{n=1}^{\infty} H_n^- \frac{x^{n+1}}{(n+1)^2} = \operatorname{Li}_3\left(\frac{2x}{1+x}\right)-\operatorname{Li}_3\left(\frac{x}{1+x}\right)-\operatorname{Li}_3\left(\frac{1+x}{2}\right)+\operatorname{Li}_3\left(\frac{1}{2}\right)-\operatorname{Li}_3(x) \qquad (20)$$

$$+\log(1+x)\left[\operatorname{Li}_2(x)+\operatorname{Li}_2\left(\frac{1}{2}\right)+\frac{1}{2}\log 2\log(1+x)\right].$$

The RHS here includes a very interesting combination of trilogarithms. We shall comment on this function in Section 3.

*Proof.* Dividing by $x$ in (7) and integrating we find

$$\sum_{n=1}^{\infty} H_n^- \frac{x^{n+1}}{(n+1)^2} = \int_0^x \left[\operatorname{Li}_2\left(\frac{1-t}{2}\right)-\operatorname{Li}_2\left(\frac{1}{2}\right)\right]\frac{dt}{t} - \log 2\int_0^x \log(1-t)\frac{dt}{t}.$$

Integrating by parts in the first integral yields (21)

$$\sum_{n=1}^{\infty} H_n^- \frac{x^{n+1}}{(n+1)^2} = \log x\left[\operatorname{Li}_2\left(\frac{1-x}{2}\right)-\operatorname{Li}_2\left(\frac{1}{2}\right)\right]+\log 2\operatorname{Li}_2(x)-\int_0^x \log\frac{1+t}{2}\log t\, \frac{dt}{1-t} \ .$$

The integral on the RHS we transform by the substitution $1-t=2y$ to get

$$-\int_0^x \log\frac{1+t}{2}\log t\, \frac{dt}{1-t} = -\int_{\frac{1}{2}}^{\frac{1-x}{2}} \log(1-2y)\log(1-y)\, \frac{dy}{y} \ .$$

For this integral we use the antiderivative given in [11, (A.3.5(1))] to write

$$-\int_0^x \log\frac{1+t}{2} \log t \, \frac{dt}{1-t} =$$

$$= \text{Li}_3\left(\frac{2x}{1+x}\right) - \text{Li}_3\left(\frac{x}{1+x}\right) - \text{Li}_3\left(\frac{1+x}{2}\right) + \text{Li}_3\left(\frac{1}{2}\right) - \text{Li}_3(x) - \log x \left[\text{Li}_2\left(\frac{1-x}{2}\right) - \text{Li}_2\left(\frac{1}{2}\right)\right]$$

$$+ \log\left(\frac{1+x}{2}\right)\left[\text{Li}_2(x) + \text{Li}_2\left(\frac{1}{2}\right) + \log^2 2\right] + \log 2 \, \text{Li}_2\left(\frac{1}{2}\right) + \frac{1}{2}\log^3 2 + \frac{1}{2}\log(2)\log^2\left(\frac{1+x}{2}\right).$$

The last line here simplifies to

$$\log(1+x)\left[\text{Li}_2(x) + \text{Li}_2\left(\frac{1}{2}\right) + \frac{1}{2}\log 2 \log(1+x)\right] - \log 2 \, \text{Li}_2(x),$$

and adding the value of this integral to (21) yields the desired result. □

## 3. Dilogarithmic and trilogarithmic relations

The above technique can be used to find more general expansions and identities for the dilogarithm and trilogarithm.

First, using Lemma 1 with $\lambda = -1$ we can write

$$\frac{-\log(1-\mu x)}{1+x} = \sum_{n=1}^{\infty}\left\{\sum_{k=1}^{n}\frac{(-\mu)^k}{k}\right\}(-1)^n x^n,$$

which turns into (2) when $\mu = -1$ and $x = -t$.

Next, notice that

$$\frac{d}{dx}\text{Li}_2\left(\frac{\mu(1+x)}{1+\mu}\right) = \frac{-\log(1-\mu x)}{1+x} + \frac{\log(1+\mu)}{1+x},$$

and therefore, integrating and adjusting the constant of integration we find

(22)

$$\sum_{n=1}^{\infty}\left\{\sum_{k=1}^{n}\frac{(-1)^{k-1}\mu^k}{k}\right\}\frac{(-1)^n x^{n+1}}{n+1} = \text{Li}_2\left(\frac{\mu(1+x)}{1+\mu}\right) - \text{Li}_2\left(\frac{\mu}{1+\mu}\right) - \log(1+\mu)\log(1+x)$$

($|\mu|<1$, $|x|<1$) . This is exactly Lewin's formula A.2.8.(1) in [11].

Notice also that

$$\frac{d}{dx}\left[\operatorname{Li}_2\left(\frac{(1+\mu)x}{1+x}\right) - \operatorname{Li}_2\left(\frac{x}{1+x}\right)\right] = \frac{-\log(1-\mu x)}{x(1+x)} = \sum_{n=1}^{\infty}\left\{\sum_{k=1}^{n}\frac{(-\mu)^k}{k}\right\}(-1)^n x^{n-1} . \quad (23)$$

Integrating this we arrive at the following:

**Proposition 10**. *We have the expansion ( $|\mu|<1$, $|x|<1$ )*

$$\operatorname{Li}_2\left(\frac{(1+\mu)x}{1+x}\right) - \operatorname{Li}_2\left(\frac{x}{1+x}\right) = \sum_{n=1}^{\infty}\left\{\sum_{k=1}^{n}\frac{(-\mu)^k}{k}\right\}\frac{(-x)^n}{n} . \quad (24)$$

Now about a property of the dilogarithm. Equation (23) says also that

$$\frac{d}{dx}\left[\operatorname{Li}_2\left(\frac{(1+\mu)x}{1+x}\right) - \operatorname{Li}_2\left(\frac{x}{1+x}\right)\right]$$

$$= \frac{-\log(1-\mu x)}{x} + \frac{\log(1-\mu x)}{1+x} = \frac{d}{dx}\operatorname{Li}_2(\mu x) + \frac{\log(1-\mu x)}{1+x} .$$

Integrating this and using equation (22) we get as an added bonus the Abel type identity for the dilogarithm (cf. [12, p.26]):

$$\operatorname{Li}_2\left(\frac{(1+\mu)x}{1+x}\right) + \operatorname{Li}_2\left(\frac{\mu(1+x)}{1+\mu}\right) - \operatorname{Li}_2(\mu x) \quad (25)$$

$$= \operatorname{Li}_2\left(\frac{x}{1+x}\right) + \operatorname{Li}_2\left(\frac{\mu}{1+\mu}\right) + \log(1+\mu)\log(1+x)$$

When $\mu = 1$ this becomes

$$\operatorname{Li}_2\left(\frac{2x}{1+x}\right) + \operatorname{Li}_2\left(\frac{1+x}{2}\right) - \operatorname{Li}_2(x) = \operatorname{Li}_2\left(\frac{x}{1+x}\right) + \operatorname{Li}_2\left(\frac{1}{2}\right) + \log 2 \log(1+x) ,$$

and with the help of Landen's identity [11]:

$$\operatorname{Li}_2\left(\frac{x}{1+x}\right) = -\frac{1}{2}\log^2(1+x) - \operatorname{Li}_2(-x)$$

we construct another one, namely,

$$\text{Li}_2\left(\frac{2x}{1+x}\right) = -\text{Li}_2\left(\frac{1+x}{2}\right) + \text{Li}_2\left(\frac{1}{2}\right) + \text{Li}_2(x) - \text{Li}_2(-x) + \log 2 \log(1+x) - \frac{1}{2}\log^2(1+x). \quad (26)$$

Among other things, (26) shows that the function $\text{Li}_2\left(\frac{2x}{1+x}\right)$ extends on the closed unit disk $|x| \leq 1$ except for $x = -1$. This function appears implicitly in the works of Ramanujan - see [1, Entry 8, p.249]. Ramanujan proved that if we define

$$f(t) = \sum_{n=1}^{\infty}\left(1 + \frac{1}{3} + \ldots + \frac{1}{2n-1}\right)\frac{t^{2n-1}}{2n-1},$$

then

$$f\left(\frac{t}{2-t}\right) = \frac{1}{8}\log^2(1-t) + \text{Li}_2(t).$$

With the substitution $x = \frac{t}{2-t}$ we can write Ramanujan's result in the form

$$\text{Li}_2\left(\frac{2x}{1+x}\right) + \frac{1}{4}\log^2\frac{1-x}{1+x} = 2\sum_{n=1}^{\infty}\left(1 + \frac{1}{3} + \ldots + \frac{1}{2n-1}\right)\frac{x^{2n-1}}{2n-1}. \quad (27)$$

We can also make similar remarks for the trilogarithm. First, we define for brevity the polynomials

$$H_n^-(\mu) = \sum_{k=1}^{n} \frac{(-\mu)^{k-1}}{k}, \quad \text{with} \quad H_n^-(1) = H_n^-.$$

Now it is easy to check that for every integer $m \geq 2$,

$$\frac{d}{dx}\left[\text{Li}_{m+1}\left(\frac{(1+\mu)x}{1+x}\right) - \text{Li}_{m+1}\left(\frac{x}{1+x}\right)\right] = \frac{1}{x(1+x)}\left[\text{Li}_m\left(\frac{(1+\mu)x}{1+x}\right) - \text{Li}_m\left(\frac{x}{1+x}\right)\right].$$

When $m = 2$ we find according to Proposition 10 and Lemma 1:

$$\frac{d}{dx}\left[\text{Li}_3\left(\frac{(1+\mu)x}{1+x}\right) - \text{Li}_3\left(\frac{x}{1+x}\right)\right] = \frac{1}{x(1+x)}\left[\text{Li}_2\left(\frac{(1+\mu)x}{1+x}\right) - \text{Li}_2\left(\frac{x}{1+x}\right)\right]$$

$$= \frac{1}{x(1+x)}\sum_{n=1}^{\infty}\frac{(-1)^{n-1}H_n^-(\mu)}{n}x^n = \frac{1}{x}\sum_{n=1}^{\infty}\left\{\sum_{k=1}^{n}\frac{H_k^-(\mu)}{k}\right\}(-1)^{n-1}x^n$$

$$= \sum_{n=1}^{\infty}\left\{\sum_{k=1}^{n}\frac{H_k^-(\mu)}{k}\right\}(-1)^{n-1}x^{n-1}.$$

Thus we come (after integration) to the remarkable expansion, analogous to (24):

**Proposition 11.** *For* $|\mu|<1, |x|<1$,

$$\text{Li}_3\left(\frac{(1+\mu)x}{1+x}\right) - \text{Li}_3\left(\frac{x}{1+x}\right) = \sum_{n=1}^{\infty}\left\{\sum_{k=1}^{n}\frac{H_k^-(\mu)}{k}\right\}\frac{(-1)^{n-1}}{n}x^n \qquad (28)$$

(*considering, if needed, a holomorphic extension of the LHS*).

Another observation: when $\mu=1$ we can write in view of equation (20)

$$\text{Li}_3\left(\frac{2x}{1+x}\right) - \text{Li}_3\left(\frac{x}{1+x}\right) - \log(1+x)\left[\text{Li}_2(x) + \text{Li}_2\left(\frac{1}{2}\right) + \frac{1}{2}\log 2\log(1+x)\right]$$

$$= \sum_{n=1}^{\infty} H_n^- \frac{x^{n+1}}{(n+1)^2} + \text{Li}_3\left(\frac{1+x}{2}\right) - \text{Li}_3\left(\frac{1}{2}\right) + \text{Li}_3(x),$$

and since the RHS is a function well defined on the closed unit disk, the function on the LHS extends to the closed unit disk.

## 4. Double integrals

Constructing a new proof that $\zeta(3)$ is irrational, F. Beukers wrote $\zeta(3)$ as a special double integral [2]. Guillera and Sondow [9] used similar double integrals to represent several classical constants and some generating functions. Here we list four more double integrals.

The generating functions considered in Propositions 5 and 8 can be represented by double integrals. Starting from

$$\frac{1}{n+k} = \int_0^1 x^{n+k-1}dx,$$

we write

$$\log 2 - H_n^- = (-1)^n \sum_{k=1}^{\infty} \frac{(-1)^{k-1}}{n+k} = (-1)^n \int_0^1 x^n \left\{\sum_{k=1}^{\infty}(-1)^{k-1}x^{k-1}\right\}dx = (-1)^n \int_0^1 \frac{x^n}{1+x}dx,$$

and for $|z|<1$,

$$\sum_{n=0}^{\infty}\left(H_n^- - \log 2\right)^2 z^n = \sum_{n=0}^{\infty}\left(\int_0^1 \frac{x^n dx}{1+x}\right)^2 z^n = \sum_{n=0}^{\infty}\left(\int_0^1 \frac{x^n dx}{1+x}\right)\left(\int_0^1 \frac{y^n dy}{1+y}\right)z^n = \sum_{n=0}^{\infty}\int_0^1\int_0^1 \frac{x^n y^n z^n dx dy}{(1+x)(1+y)}$$

$$= \int_0^1\int_0^1 \left\{\sum_{n=0}^{\infty}(xyz)^n\right\} \frac{dx\,dy}{(1+x)(1+y)} = \int_0^1\int_0^1 \frac{dx\,dy}{(1-xyz)(1+x)(1+y)}.$$

Thus from Proposition 5 we have the following representation.

**Corollary 12.** *For* $|z| \leq 1$ (29)

$$g(z) = \int_0^1\int_0^1 \frac{dx\,dy}{(1-xyz)(1+x)(1+y)} = \frac{1}{1-z}\left\{\text{Li}_2(z) + \log^2 2 - 2\left[\text{Li}_2\left(\frac{1+z}{2}\right) - \text{Li}_2\left(\frac{1}{2}\right)\right]\right\}$$

*with* $g(1) = \log 2$ *and* $g(-1) = \dfrac{\pi^2}{24}$.

*Further, integrating for* $z$ *(i.e. using (17))*,

$$G(z) \equiv -\int_0^1\int_0^1 \frac{\log(1-xyz)\,dx\,dy}{xy(1+x)(1+y)} = \int_0^z \frac{\text{Li}_2(t)}{1-t}dt + \log(1-z)\left[2\text{Li}_2\left(\frac{1+z}{2}\right) - \frac{\pi^2}{6}\right] \quad (30)$$

$$+ 2\log(1+z)\left(\log^2(1-z) - \log^2 2\right) + 2\log 2\left[\text{Li}_2\left(\frac{1+z}{2}\right) - \text{Li}_2\left(\frac{1}{2}\right)\right]$$

$$- 2\log 2 \log^2(1-z) + 4\log(1-z)\text{Li}_2\left(\frac{1-z}{2}\right) - 4\left[\text{Li}_3\left(\frac{1-z}{2}\right) - \text{Li}_3\left(\frac{1}{2}\right)\right],$$

*and the values* $G(\pm 1)$ *are given above in Proposition 8.*

Two similar double integrals were recently evaluated by this author. In the solution of Problem H-691 in the Fibonacci Quarterly [4] we found

$$\int_0^1\int_0^1 \frac{x^2 y^2 dx\,dy}{(1+x^2 y^2)(1+x)(1+y)} = \frac{7}{8}(\ln 2)^2 + \frac{\pi}{8}\ln 2 - \frac{G}{2} - \frac{\pi^2}{48}, \quad (31)$$

where $G$ is the Catalan constant.

Another double integral appeared in [3],

$$\int_0^1\int_0^1 \frac{\ln(1+xy)\,dx\,dy}{(1+x)(1+y)} = \frac{\pi^2}{12}\log 2 + \frac{1}{3}\log^3 2 - \frac{1}{2}\zeta(3). \quad (32)$$

## 5. Work done by De Doelder

In the second part of his article [5] De Doelder studied the series

$$d_1(x) = \sum_{n=1}^{\infty}\left[\psi\left(\frac{n+1}{2}\right) - \psi\left(\frac{n}{2}\right)\right]\frac{x^n}{n} \quad \text{and}$$

$$d_2(x) = \sum_{n=1}^{\infty}\left[\psi\left(\frac{n+1}{2}\right) - \psi\left(\frac{n}{2}\right)\right]^2 \frac{x^n}{n},$$

where the second series was evaluated only for $x = 1$ and $x = -1$. In view of equation (1) we see that the series $d_1(x)$ is equivalent to the series (8) and the evaluations of $d_2(\pm 1)$ are similar to (18) and (19).